\documentstyle[12pt]{article}

\newtheorem{t.}{Theorem}[section]
\newtheorem{d.}[t.]{Definition}
\newtheorem{l.}[t.]{Lemma}
\newtheorem{p.}[t.]{Proposition}
\newtheorem{c.}[t.]{Corollary}
\newtheorem{e.}[t.]{Example}
\newtheorem{r.}[t.]{Remark}
\begin{document}
\title {More on Reverse Triangle Inequality in Inner Product Spaces}
\author{A. H. Ansari and M. S. Moslehian}
\date{}
\maketitle

\begin{abstract}
Refining some results of S. S. Dragomir, several new reverses of the generalized triangle inequality in inner product spaces are given. Among several results, we establish some reverses for the Schwarz inequality. In particular, it is proved that if $a$ is a unit vector in a real or complex inner product space $(H;\langle .,.\rangle)$, $r, s>0, p\in(0,s], D=\{x\in H,\|rx-sa\|\leq p\}, x_1, x_2\in D-\{0\}$ and $ \alpha_{r,s}=\min\{\frac{r^2\|x_k\|^2-p^2+s^2}{2rs\|x_k\|}: 1\leq k\leq 2 \}$, then
$$\frac{\|x_1\|\|x_2\|-Re\langle x_1,x_2\rangle}{(\|x_1\|+\|x_2\|)^2}\leq \alpha_{r,s}.$$
\end{abstract}

\section{Introduction.}
	
It is interesting to know under which conditions the triangle inequality went the other way in a normed space $X$; in other words, we would like to know if there is a constant $c$ with the property that $c\sum_{k=1}^n\|x_k\|\leq\|\sum_{k=1}^n x_k\|$ for any finite set $x_1,\cdots,x_n\in X$. M. Nakai and T. Tada $\cite{N-T}$ proved that the normed spaces with this property are precisely those of finite dimensional.

The first authors investigating reverse of the triangle inequality in inner product spaces were J. B. Diaz and F. T. Metcalf $\cite{D-M}$ by establishing the following result as an extension of an inequality given by M. Petrovich $\cite{PET}$ for complex numbers:

{\bf Diaz-Metcalf Theorem.} Let $a$ be a unit vector in an inner product space $(H;\langle .,.\rangle)$. Suppose the vectors $x_k \in H, k\in\{1,\cdots,n\}$ satisfy
\begin{eqnarray*}
0\leq r\leq\frac{Re\langle x_k,a\rangle}{\|x_k\|}, k\in\{1,\cdots,n\}
\end{eqnarray*}
Then 
\begin{eqnarray*}
r\sum_{k=1}^n\|x_k\|\leq\|\sum_{k=1}^n x_k\|.
\end{eqnarray*}
where equality holds if and only if
\begin{eqnarray*}
\sum_{k=1}^n x_k =r\sum_{k=1}^n \|x_k\|a.
\end{eqnarray*}

Inequalities related to the triangle inequality are of special interest; cf. Chapter XVII of $\cite{M-P-F}$. They may be applied to get interesting inequalities in complex numbers or to study vector-valued integral inequalities $\cite{DRA1}, \cite{DRA2}$.

Using several ideas and following the terminology of $\cite{DRA1}$ and $\cite{DRA2}$ we modify or refine some results of Dragomir and ours $\cite{A-M}$ and get some new reverses of triangle inequality. Among several results, we show that if $a$ is a unit vector in a real or complex inner product space $(H;\langle .,.\rangle)$, $x_k \in H-\{0\}, 1\leq k\leq n$, $\alpha=\min\{\|x_k\|: 1\leq k\leq n\}, p\in(0,\sqrt{\alpha^2+1}), \max\{\|x_k-a\|: 1\leq k\leq n\}\leq p$ and $\beta=\min\{\frac{\|x_k\|^2-p^2+1}{2\|x_k\|}: 1\leq k\leq n\}$, then 
$$\sum_{k=1}^n\|x_k\| -\|\sum_{k=1}^n x_k\|\leq\frac{1-\beta}{\beta}Re\langle\sum_{k=1}^n x_k ,a\rangle.$$ 
We also examin some reverses for the celebrated Schwarz inequality. In particular, it is proved that if $a$ is a unit vector in a real or complex inner product space $(H;\langle .,.\rangle)$, $r, s>0, p\in(0,s], D=\{x\in H,\|rx-sa\|\leq p\}, x_1, x_2\in D-\{0\}$ and $ \alpha_{r,s}=\min\{\frac{r^2\|x_k\|^2-p^2+s^2}{2rs\|x_k\|}: 1\leq k\leq 2 \}$, then
$$\frac{\|x_1\|\|x_2\|-Re\langle x_1,x_2\rangle}{(\|x_1\|+\|x_2\|)^2}\leq \alpha_{r,s}$$

Throughout the paper $(H;\langle .,.\rangle)$ denotes a real or complex inner product space. We use repeatedly the Cauchy-Schwarz inequality without mentioning it. The reader is referred to $\cite{RAS}, \cite{DRA3}$ for the terminology on inner product spaces.

\section{Reverse of Triangle inequality.}

We start this section by pointing out the following theorem of $\cite{A-M}$ which is a modification of theorem 3 of $\cite{DRA2}$.

\begin{t.} Let $a_{1},\ldots ,a_{m}$ be orthonormal vectors in the complex inner product space $(H;\langle .,.\rangle)$. Suppose that for $1\leq t\leq m ,r_t,\rho_t\in R$ and that the vectors $x_k \in H, k\in\{1,\cdots,n\}$ satisfy 
\begin{equation}
0\leq r_t ^2 \|x_k\|\leq Re\langle x_k,r_t a_t\rangle, 0\leq \rho_t ^2 \|x_k\| \leq Im\langle x_k,\rho_t a_t \rangle, 1\leq t\leq m
\end{equation}
Then 
\begin{equation}
(\sum_{t=1}^m (r_t ^2+\rho_t ^2))^\frac{1}{2}\sum_{k=1}^n\|x_k\|\leq\|\sum_{k=1}^n x_k\|
\end{equation}
and the equality holds in $(2)$ if and only if
\begin{equation}
\sum_{k=1}^n x_k =\sum_{k=1}^n \|x_k\|\sum_{t=1}^m (r_t+i\rho_t)a_t.
\end{equation}
\end{t.}

The following theorem is a strengthen of Corollary 1 of $\cite{DRA2}$ and a generalization of Theorem 2 of $\cite{A-M}$.

\begin{t.} Let $a$ be a unit vector in the complex inner product space $(H;\langle .,.\rangle)$. Suppose that the vectors $x_k \in H-\{0\}, k\in\{1,\cdots,n\}$ satisfy 
$$\max\{\|rx_k-sa\|: 1\leq k\leq n\}\leq p, \max\{\|r'x_k-is'a\|: 1\leq k\leq n\}\leq q$$ where $r,r',s,s'>0$ and 
$$p\leq((r\alpha)^2+s^2)^\frac{1}{2}, q\leq((r'\alpha)^2+s'^2)^\frac{1}{2}$$
and $$\alpha=\min\{\|x_k\|: 1\leq k\leq n\}.$$
Let
$$\alpha_{r,s}=\min\{\frac{r^2\|x_k\|^2-p^2+s^2}{2rs\|x_k\|}: 1\leq k\leq n\},$$
$$\beta_{r',s'}=\min\{\frac{r'^2\|x_k\|^2-q^2+s'^2}{2r's'\|x_k\|}: 1\leq k\leq n\}$$
Then
$$(\alpha_{r,s}^2+\beta_{r',s'}^2)^\frac{1}{2}\sum_{k=1}^n\|x_k\|\leq\|\sum_{k=1}^n x_k\|$$
and the equality holds if and only if
$$ \sum_{k=1}^n x_k= (\alpha_{r,s}+i\beta_{r',s'})\sum_{k=1}^n\|x_k\|a.$$
\end{t.}

{\bf Proof.} From the first inequality above we infer that 
$$\langle rx_k-sa,rx_k-sa\rangle\leq p^2$$
$$r^2\|x_k\|^2+s^2-p^2\leq2Re\langle rx_k,sa\rangle.$$ Then
$$\frac{r^2\|x_k\|^2+s^2-p^2}{2rs\|x_k\|}\|x_k\|\leq Re\langle x_k,a\rangle.$$Similarly
$$\frac{r'^2 \|x_k\|^2-q^2+s'^2}{2r's' \|x_k\|} \|x_k\|\leq Im\langle x_k,a\rangle$$ consequently
$$\alpha_{r,s}\|x_k\|\leq Re\langle x_k,a\rangle$$
and
$$\beta_{r',s'}\|x_k\|\leq Im\langle x_k,a\rangle.$$
Applying Theorem 2.1 for $m=1, r_1=\alpha_{r,s}$ and $\rho_1=\beta_{r',s'}$ we deduce desired inequality. $\Box$

The next result is an extension of Corollary 3 of $\cite{A-M}$.

\begin{c.} Let $a$ be a unit vector in the complex inner product space $(H;\langle .,.\rangle)$. Suppose that $x_k \in H, k\in\{1,\cdots,n\}, \max\{\|rx_k-sa\|: 1\leq k\leq n\}\leq r, \max\{\|rx_k-isa\|: 1\leq k\leq n\}\leq s$ where $r>0, s>0$ and $\alpha=\min\{\|x_k\|: 1\leq k\leq n\}.$
Then 
$$\frac{r\alpha}{s\sqrt{2}}\sum_{k=1}^n\|x_k\|\leq\|\sum_{k=1}^n x_k\|.$$ The equality holds if and only if
$$\sum_{k=1}^n x_k= r\alpha\frac{(1+i)}{2s}\sum_{k=1}^n\|x_k\|a.$$
\end{c.}

{\bf Proof.} Apply Theorem 2.2 with $r=r', s=s', p=r, q=s$. Note that $\alpha_{r,s}=\frac{r\alpha}{2s}=\beta_{r',s'}. \Box$ 

\begin{t.} Let $a$ be a unit vector in $H;\langle .,.\rangle)$. Suppose that the vectors $x_k \in H, k\in\{1,\cdots,n\}$ satisfy 
$$\max\{\|rx_k-sa\|: 1\leq k\leq n\}\leq p<((r\alpha)^2+s^2)^\frac{1}{2}$$ 
where $r>0, s>0$ and
$$\alpha=\displaystyle{\min_{1\leq k\leq n}}\|x_k\|.$$
Let
$$\alpha_{r,s}=\min\{\frac{r^2\|x_k\|^2-p^2+s^2}{2rs\|x_k\|}:1\leq k\leq n \}.$$
Then 
$$\alpha_{r,s}\sum_{k=1}^n\|x_k\|\leq\|\sum_{k=1}^n x_k\|.$$
Moreover, the equality holds if and only if
$$\sum_{k=1}^n x_k=\alpha_{r,s}\sum_{k=1}^n\|x_k\|a.$$
\end{t.}

{\bf Proof}. Proof is similar to that of Theorem 2.2 in which we use Theorem 2.1 with $m=1, \rho_1=0. \Box$ 

\begin{t.} Let $a$ be a unit vector in $(H;\langle .,.\rangle)$. Suppose that $r, s>0$, and vectors $x_k \in H-\{0\}, k\in\{1,\cdots,n\}$ satisfy 
$$\sum_{k=1}^n x_k=0.$$
Then 
$$\sqrt{r^2\alpha^2+s^2}\leq \max\{\|rx_k-sa\|:1\leq n\}$$
where
$$\alpha=\displaystyle{\min_{1\leq k\leq n}}\|x_k\|.$$
\end{t.}

{\bf Proof}. Let $ p=\max\{\|rx_k-sa\|:1\leq k\leq n\}$. If $p<\sqrt{r^2\alpha^2+s^2}$, then using Theorem 2.4 we get
$$\alpha_{r,s}\sum_{k=1}^n\|x_k\|\leq\|\sum_{k=1}^n x_k\|=0.$$
Hence $\alpha_{r,s}=0$. On the other hand $\frac{p^2-s^2}{r^2}<\alpha^2$, so
$$\alpha_{r,s}=\min\{\frac{r^2\|x_k\|^2-p^2+s^2}{2rs\|x_k\|}:1\leq k\leq n \}>0$$
a contradiction. $\Box$

\begin{t.} Let $a_{1},\ldots ,a_m$ be orthonormal vectors in the complex inner product space $(H;\langle .,.\rangle),M_t\geq m_t > 0, L_t\geq \ell_t>0, 1\leq t\leq m$ and $x_k \in H-\{0\}, k\in\{1,\cdots,n\}$ such that
$$Re\langle M_ta_t-x_k,x_k-m_ta_t\rangle\geq 0, Re\langle L_tia_t-x_k,x_k-\ell_tia_t\rangle\geq0$$
or equivalently
$$\|x_k-\frac{m_t+M_t}{2}a_t\|\leq\frac{M_t-m_t}{2},\|x_k-\frac{L_t+\ell_t}{2}ia_t\|\leq\frac{L_t-\ell_t}{2}$$
for all $1\leq k\leq n$ and $1\leq t\leq m$. Let
$$\alpha_{m_t,M_t}=\min\{\frac{\|x_k\|^2+m_tM_t}{(m_t+M_t)\|x_k\|}: 1\leq k\leq n\}, 1\leq t\leq m$$
and
$$\alpha_{\ell_t,L_t}=\min\{\frac{\|x_k\|^2+\ell_tL_t}{(m_t+M_t)\|x_k\|}: 1\leq k\leq n\}, 1\leq t\leq m$$ 
then 
$$(\sum_{t=1}^m \alpha_{m_t,M_t}^2+\alpha_{\ell_t,L_t}^2)^\frac{1}{2}\sum_{k=1}^n\|x_k\|\leq\|\sum_{k=1}^n x_k\|.$$
The equality holds if and only if
$$\sum_{k=1}^n x_k=(\sum_{k=1}^n\|x_k\|)\sum_{t=1}^m(\alpha_{m_t,M_t}+i\alpha_{\ell_t,L_t})a_t.$$
\end{t.}

{\bf Proof.} Given $1\leq t\leq m$ and all $1\leq k\leq n$, it follows from $\|x_k-\frac{m_t+M_t}{2}a_t\|\leq\frac{M_t-m_t}{2}$ that
$$\|x_k\|^2+m_tM_t\leq(m_t+M_t)Re\langle x_k, a_t\rangle.$$
Then
$$\frac{\|x_k\|^2+m_tM_t}{(m_t+M_t)\|x_k\|}\|x_k\|\leq Re\langle x_k, a_t\rangle$$
and so 
$$\alpha_{m_t,M_t}\|x_k\|\leq Re\langle x_k, a_t\rangle.$$
Similarly from the second inequality we deduce that
$$\alpha_{\ell_t,L_t}\|x_k\|\leq Im\langle x_k, a_t\rangle.$$
Applying Theorem 2.5 for $ r_t=\alpha_{m_t,M_t}$ and $\rho_t=\alpha_{\ell_t,L_t}$, we obtain the required inequality. $\Box$ 

We will need theorem 7 of $\cite{DRA1}$. We mention it for the sake of completeness.

\begin{t.} Let $a$ be a unit vector in $(H;\langle .,.\rangle)$, and $x_k \in H-\{0\}, k\in\{1,\cdots,n\}$. If $r_k\geq 0, k\in\{1,\cdots,n\}$ such that
$$\|x_k\|-Re\langle x_k,a\rangle\leq r_k$$
then 
$$\sum_{k=1}^n\|x_k\| -\|\sum_{k=1}^n x_k\|\leq \sum_{k=1}^n r_k.$$
The equality holds if and only if 
$$\sum_{k=1}^n\|x_k\|\geq \sum_{k=1}^n r_k $$ and 
$$ \sum_{k=1}^n x_k=(\sum_{k=1}^n\|x_k\|-\sum_{k=1}^n r_k)a.$$
\end{t.}

\begin{t.} Let $a$ be a unit vector in $(H;\langle .,.\rangle)$ and $x_k \in H-\{0\}, k\in\{1,\cdots,n\}.$
Let
$$\alpha=\min\{\|x_k\|: 1\leq k\leq n\}, p\in(0,\sqrt{\alpha^2+1}), \max\{\|x_k-a\|:1\leq k\leq n\}\leq p$$
and 
$$\beta=\min\{\frac{\|x_k\|^2-p^2+1}{2\|x_k\|}:1\leq k\leq n\}.$$
Then
$$\sum_{k=1}^n\|x_k\| -\|\sum_{k=1}^n x_k\|\leq\frac{1-\beta}{\beta}Re\langle\sum_{k=1}^n x_k ,a\rangle.$$
The equality holds if and only if 
$$\sum_{k=1}^n\|x_k\|\geq\frac{1-\beta}{\beta} Re\langle\sum_{k=1}^n x_k ,a\rangle$$
and $$\sum_{k=1}^n x_k=(\sum_{k=1}^n\|x_k\|-\frac{1-\beta}{\beta}Re\langle\sum_{k=1}^n x_k ,a\rangle)a.$$
\end{t.}

{\bf Proof.} Since $\max\{\|x_k-a\|:1\leq k\leq n\}\leq p$, we have
$$\langle x_k-a,x_k-a\rangle\leq p^2$$
$$\| x_k\|^2+1-p^2\leq2Re\langle x_k, a\rangle$$
$$\frac{\|x_k\|^2-p^2+1}{2\|x_k\|}\|x_k\|\leq Re\langle x_k,a\rangle$$
$$\beta\| x_k\|\leq Re\langle x_k ,a\rangle$$
$$\| x_k\|\leq\frac{1}{\beta}Re\langle x_k ,a\rangle$$
for all $k\in\{1,\cdots,n\}$. Then
$$\| x_k\|-Re\langle x_k ,a\rangle\leq\frac{1-\beta}{\beta}Re\langle x_k ,a\rangle, k\in\{1,\cdots,n\}.$$
Applying Theorem 2.7 for $r_k=\frac{1-\beta}{\beta}Re\langle x_k ,a\rangle, k\in\{1,\cdots,n\}$, we deduce the desired inequality. $\Box$

As a corollary, we obtain a result similar to Theorem 9 of $\cite{DRA1}$:

\begin{c.} Let $a$ be a unit vector in $(H;\langle .,.\rangle)$ and $x_k \in H-\{0\}, k\in\{1,\cdots,n\}$.
Let 
$$\max\{\|x_k-a\|: 1\leq k\leq n\}\leq 1$$
and
$$\alpha=\min\{\|x_k\|: 1\leq k\leq n\}.$$
Then 
$$\sum_{k=1}^n\|x_k\| -\|\sum_{k=1}^n x_k\|\leq\frac{2-\alpha}{\alpha}Re\langle\sum_{k=1}^n x_k ,a\rangle.$$
The equality holds if and only if
$$\sum_{k=1}^n\|x_k\|\geq\frac{2-\alpha}{\alpha}Re\langle\sum_{k=1}^n x_k ,a\rangle$$
and $$\sum_{k=1}^n x_k=(\sum_{k=1}^n\|x_k\|-\frac{2-\alpha}{\alpha}Re\langle\sum_{k=1}^n x_k ,a\rangle)a.$$
\end{c.}

{\bf Proof.} Apply Theorem 2.8 with $\beta=\frac{\alpha}{2}. \Box$ 

\begin{t.} Let $a$ be a unit vector in $(H;\langle .,.\rangle)$, $ M\geq m>0$ and $x_k \in H-\{0\}, k\in\{1,\cdots,n\}$ such that 
$$Re\langle Ma-x_k,x_k-ma\rangle\geq 0$$
or equivalently
$$\|x_k-\frac{m+M}{2}a\|\leq\frac{M-m}{2}.$$
Let
$$\alpha_{m,M}=\min\{\frac{\|x_k\|^2+mM}{(m+M)\|x_k\|}:1\leq k\leq n\}.$$ 
Then 
$$\sum_{k=1}^n\|x_k\| -\|\sum_{k=1}^n x_k\|\leq\frac{1-\alpha_{m,M}}{\alpha_{m,M}}Re\langle\sum_{k=1}^n x_k ,a\rangle.$$
The equality holds if and only if 
$$\sum_{k=1}^n\|x_k\|\geq\frac{1-\alpha_{m,M}}{\alpha_{m,M}} Re\langle\sum_{k=1}^n x_k ,a\rangle$$
and 
$$\sum_{k=1}^n x_k=(\sum_{k=1}^n\|x_k\|-\frac{1-\alpha_{m,M}}{\alpha_{m,M}}Re\langle\sum_{k=1}^n x_k ,a\rangle)a.$$
\end{t.}

{\bf Proof.} For each $1\leq k\leq n$, it follows from the inequality
$$\|x_k-\frac{m+M}{2}a\|\leq\frac{M-m}{2}$$
that 
$$\langle x_k-\frac{m+M}{2}a,x_k-\frac{m+M}{2}a\rangle\leq(\frac{M-m}{2})^2.$$
Hence
$$\|x_k\|^2+mM\leq(m+M)Re\langle x_k, a\rangle.$$
So that
$$\alpha_{m,M}\|x_k\|\leq Re\langle x_k, a\rangle$$
consequently
$$\| x_k\|-Re\langle x_k ,a\rangle\leq\frac{1-\alpha_{m,M}}{\alpha_{m,M}}Re\langle x_k ,a\rangle.$$
Now apply Theorem 2.7 for $r_k=\frac{1-\alpha_{m,M}}{\alpha_{m,M}}Re\langle x_k ,a\rangle, k\in\{1,\cdots,n\}. \Box$

\section{Reverses of Schwarz inequality.}

In this section we provide some reverses of the Schwarz inequality. The first theorem is an extension of Proposition 5.1 of $\cite{DRA1}$.

\begin{t.} Let $a$ be a unit vector in $(H;\langle .,.\rangle)$. Suppose that $r, s>0, p\in(0,s]$ and $$D=\{x\in H,\|rx-sa\|\leq p\}.$$
If $0\neq x_1\in D, 0\neq x_2\in D$, then 
$$\frac{\|x_1\|\|x_2\|-Re\langle x_1,x_2\rangle}{(\|x_1\|+\|x_2\|)^2}\leq\frac{1}{2}(1-(\frac{r^2\|x_1\|^2-p^2+s^2}{2rs\|x_1\|})^2)$$ or
$$\frac{\|x_1\|\|x_2\|-Re\langle x_1,x_2\rangle}{(\|x_1\|+\|x_2\|)^2}\leq\frac{1}{2}(1-(\frac{r^2\|x_2\|^2-p^2+s^2}{2rs\|x_2\|})^2)$$
\end{t.}

{\bf Proof.} Put $\alpha_{r,s}=\min\{\frac{r^2\|x_k\|^2-p^2+s^2}{2rs\|x_k\|}:1\leq k\leq 2 \}$. 
By Theorem 2.4, we obtain
$$\alpha_{r,s}(\|x_1\|+\|x_2\|)\leq\|x_1+x_2\|.$$
Then
$$\alpha_{r,s}^2(\|x_1\|^2+2\|x_1\|\|x_2\|+\|x_2\|^2)\leq\|x_1\|^2+2Re\langle x_1,x_2\rangle+\|x_2\|^2.$$
Set $\alpha_{r,s}^2=1-t^2$. Then
$$\frac{\|x_1\|\|x_2\|-Re\langle x_1,x_2\rangle}{(\|x_1\|+\|x_2\|)^2}\leq\frac{1}{2}t^2,$$
namely
$$\frac{\|x_1\|\|x_2\|-Re\langle x_1,x_2\rangle}{(\|x_1\|+\|x_2\|)^2}\leq\frac{1}{2}(1-\alpha_{r,s}^2). \Box$$

\begin{c.} Let $a$ be a unit vector in $(H;\langle .,.\rangle)$. Suppose that $r, s>0$ and
$$D=\{x\in H,\|rx-sa\|\leq s\}.$$
If $x,y\in D$ and $0<\|x\|<\|y\|$, then 
$$\frac{\|x\|\|y\|-Re\langle x,y\rangle}{(\|x\|+\|y\|)^2}\leq\frac{1}{2}(1-(\frac{r\|x\|}{2s})^2).$$
\end{c.}

{\bf Proof.} In the notation of the proof of Theorem 3.1 we get from $p=s, x_1=x, x_2=y$ that $\alpha_{r,s}=\frac{r\|x\|}{2s}$. Now apply Theorem 3.1. $\Box$

\begin{c.} Let $a$ be a unit vector in $(H;\langle .,.\rangle)$, $M\geq m>0$ and $x_k \in H-\{0\}, k=1,2$ such that
$$Re\langle Ma-x_k,x_k-ma\rangle\geq 0$$
or equivalently, $$\|x_k-\frac{m+M}{2}a\|\leq\frac{M-m}{2}.$$
Then 
$$\frac{\|x_1\|\|x_2\|-Re\langle x_1,x_2\rangle}{(\|x_1\|+\|x_2\|)^2}\leq\frac{1}{2}(1-(\frac{\|x_1\|^2+mM}{(m+M)\|x_1\|})^2)$$
or
$$\frac{\|x_1\|\|x_2\|-Re\langle x_1,x_2\rangle}{(\|x_1\|+\|x_2\|)^2}\leq\frac{1}{2}(1-(\frac{\|x_2\|^2+mM}{(m+M)\|x_2\|})^2).$$
\end{c.}

{\bf Proof.} Put $r=1, s=\frac{m+M}{2}, p=\frac{M-m}{2}, x=x_1$ and $y=x_2$ in Theorem 3.1. $\Box$ 

{\bf Acknowledgement.} The authors would like to thank the referees for their valuable suggestions.

Arsalan Hojjat Ansari\\
Mohammad Sal Moslehian\\
Dept. of Math., Ferdowsi Univ., P. O. Box 1159, Mashhad 91775, Iran\\
E-mail: moslehian@ferdowsi.um.ac.ir\\
Home: http://www.um.ac.ir/$\sim$moslehian/
\end{document}